\documentclass[12pt]{article}
\usepackage[T2A]{fontenc}
\usepackage[cp1251]{inputenc}
\usepackage{amsfonts}
\usepackage{eufrak}
\usepackage{amssymb}

\usepackage{graphicx}
\usepackage{graphics}

\usepackage{cite}
\usepackage{color}

\usepackage{epstopdf}
\usepackage{epsfig}
\usepackage{subfigure}

\begin{document}


\centerline{\textbf{RANDOM WALKS ON }}

\centerline{\textbf{DISCRETE ABELIAN GROUPS}}

\bigskip

\centerline{\textbf{Margaryta Myronyuk}}

\begin{abstract}
In the present paper we find necessary and sufficient conditions for
recurrence of random walks on arbitrary subgroups of the group of
rational numbers $\mathbb{Q}$. \end{abstract}

\bigskip


\bigskip


\bigskip
\bigskip
\bigskip

\section{Introduction}

Let $(\Omega, \mathfrak{A}, P)$ be a probability space, $X$ be a
countable discrete Abelian groups, $\mu$ be a distribution on $X$.
Recall that a random walk on a group $X$ generated by the
distribution $\mu$ is a sequence
$$ S_n=\xi_1+\dots+\xi_n, \quad n=1, 2, \dots,
$$
where $\xi_j$ are independent identically distributed with the
distribution $\mu$ random variables defined on $(\Omega,
\mathfrak{A}, P)$ with values on $X$. The random walk on the group
$X$ is said to be recurrent if all elements of  $X$ are recurrent
i.e. for every $x\in X$ the equality
\begin{equation}
\label{eq1b} P\{\omega\in\Omega: S_n(\omega)=x \mbox{\ for
infinitely many indices }n\in\mathbf{N}\}=1
\end{equation} holds.

Denote by $\mathbb{Z}$ the additive group of integers, by
$\mathbb{Q}$ the additive group of rational numbers considering in
the discrete topology, and by $\mathbb{Z}(k)$ the finite cyclic
group. R. M. Dudley (\cite{Dudley}) proved that there exists a
recurrent random walk on a countable Abelian group $X$ iff $X$
contains no subgroup isomorphic to $\mathbb{Z}^3$. Dudley's proof is
not constructive and therefore it is natural to look for other
effective recurrence conditions.

Such conditions on $\mu$ were studied

(a) on the weak direct product $\mathbb{Z}(2)^{\aleph_0*}$
(\cite{Darling-Erdos}),

(b) on the weak direct product $\mathbf{P}_{i\in
\mathbb{N}}^*\mathbb{Z}(k_i)$ (\cite{Flatto-Pitt}),

(c) on the factor group $\mathbb{Q}/\mathbb{Z}$ and its subgroups
(\cite{Flatto-Pitt}),

(d) on the weak direct product $\mathbb{Z}(k)^{\aleph_0*}$
(\cite{Ferage-Molchanov}),

(e) on subgroups $H_p=\left\{ {m\over p^n}: n=1,2,...; m\in
\mathbb{Z} \right\}$ of the group  $\mathbb{Q}$
(\cite{Ferage-Molchanov}),

(f) on groups of the form $\mathbb{Z}^k\times \mathbf{P}_{i\in
\mathbb{N}}^*\mathbb{Z}(k_i)$ (\cite{Kasym}).

In the present paper we find necessary and sufficient conditions for
recurrence of random walks on arbitrary subgroups of the group of
rational numbers $\mathbb{Q}$.

Let $X$ be a locally compact second countable Abelian group,
$Y=X^\ast$ be its character group, and $(x,y)$ be the value of a
character $y \in Y$ at an element $x \in X$. Let
$$\widehat \mu(y) = \int_X (x, y) d\mu(x)$$ be the characteristic function
of a distribution $\mu$ on $X$. Denote by $m_X$ the Haar measure on
 $X$.

The following recurrence criterion was proved in
\cite{Kesten-Spitzer}.

\bigskip

\textbf{Theorem А} (\cite{Kesten-Spitzer}). \textit{Let $X$ be a
countable discrete Abelian group. Let $\mu$ be a distribution on
$X$. A random walk defined by the distribution $\mu$ is recurrent
iff}
\begin{equation}\label{1}
    \int_Y Re{1\over 1-\hat\mu(y)} dm_Y = \infty
\end{equation}

\bigskip

In cases (b), (c), (d), (f) authors of corresponding papers use
Theorem A to obtain necessary and sufficient conditions for
recurrence of random walks. In these cases the character groups have
quite simple structure.

The case (e) is more complicated. Authors of \cite{Ferage-Molchanov}
state  that Theorem A is useless for
obtaining necessary and sufficient conditions for recurrence of
random walks on a group $X$ if its character group $Y$ has a
complicated structure as, for example, the group $H_p=\left\{
{m\over p^n}: n=1,2,...; m\in \mathbb{Z} \right\}$. In this case the
character group $Y$ is a $p$-adic solenoid. In \cite{Ferage-Molchanov}
necessary and sufficient conditions for recurrence of
random walks on a group $H_p$ were obtained without using of Theorem A.

In the present paper we do use Theorem A and find necessary and
sufficient conditions for recurrence of random walks on arbitrary
subgroups of the group of rational numbers $\mathbb{Q}$. In this
case the character group $Y$ is a $\boldmath{a}$-adic solenoid. Note
that the results of the paper \cite{Ferage-Molchanov} for the groups
$H_p$
follows directly from our paper.

\section{Notation and definitions}

In the present paper we use some results from the Pontryagin duality
theory (see \cite{HeRo1}).

In the paper we consider random walks on arbitrary subgroups of the
group of rational numbers not isomorphic to $\mathbb{Z}$.

Let $\boldmath{a}=(a_1,a_2,...)$ be a sequence of integers where all
$a_j>1$. Consider a group
\begin{equation}\label{def1}
    H_{\boldmath{a}}=\left\{ {m\over a_1a_2...a_n}: n=1,2,...; m\in
\mathbb{Z} \right\}.
\end{equation}
It is well-known that any subgroup of the group $\mathbb{Q}$ not
isomorphic to $\mathbb{Z}$ has the form (\ref{def1}) for some
$\boldmath{a}=(a_1,a_2,...)$. Particularly, if all $a_j=p$ then we
obtain the group of $H_p=\left\{ {m\over p^n}: n=1,2,...; m\in
\mathbb{Z} \right\}$.

In order to apply Theorem А, we have to describe the character group
of the group $H_{\boldmath{a}}$.

Let $\Delta_{\boldmath{a}}$ be a group of $\boldmath{a}$-adic
integers. Consider the group
$\mathbb{R}\times\Delta_{\boldmath{a}}$. Let $B$ the subgroup of
$\mathbb{R}\times\Delta_{\boldmath{a}}$ of the form
$B=\{(n,n\mathbf{u}\}_{n=-\infty}^{\infty}$, where
$\mathbf{u}=(1,0,0,...,0,...)\in \Delta_{\boldmath{a}}$.The factor
group
$\Sigma_{\boldmath{a}}=\mathbb{R}\times\Delta_{\boldmath{a}}/B$ is
called the $\boldmath{a}$-adic solenoid. The group
$\Sigma_{\boldmath{a}}$ is compact and connected (see
\cite[$\S$10]{HeRo1}). The group $\Sigma_{\boldmath{a}}$ is
topologically isomorphic to the character group of the group
$H_{\boldmath{a}}$ (see \cite[$($25.3)]{HeRo1}).

Denote by $\mathbb{T}$  the circle group (the one-dimensional
torus), i.e. $\mathbb{T}=\{z\in \mathbb{C}: \ |z|=1\}$. It is
convenient for us to use another consideration of the
${\boldmath{a}}$-adic solenoid as a subgroup of the
infinite-dimensional torus $\mathbb{T}^{\aleph_0}$.

Consider the mapping $f: \mathbb{R}\times\Delta_{\boldmath{a}}
\longrightarrow \mathbb{T}^{\aleph_0}$, defined by

$$ f(t,\mathbf{y})\longmapsto z=(z_1,z_2,...), \ z_j=\exp \left({2\pi i\over a_1\cdot\cdot\cdot a_j}
\left(t-(y_0+a_1y_1+...+a_1a_2...a_{j-1}y_{j-1}) \right) \right),$$

\noindent where $t\in\mathbb{R}$, $\mathbf{y}=(y_0,y_1,...) \in
\Delta_{\boldmath{a}}$.

It is not difficult to verify that $f$ is a continuous homomorphism,
${\rm Ker} f=B$ and $G={\rm Im} f={\{z=(z_1, z_2, \dots, z_n,
\dots)\in\mathbb{T}^{\aleph_0}: z_k^{a_k}=z_{k-1}\}}$ is a closed
subgroup of the infinite-dimensional torus $\mathbb{T}^{\aleph_0}$.
Then $G\cong\Sigma_{\boldmath{a}}$.

The consideration of the $\boldmath{a}$-adic solenoid
$\Sigma_{\boldmath{a}}$ as the subgroup
$$G={\{z=(z_1, z_2, \dots, z_n, \dots)\in\mathbb{T}^{\aleph_0}:
z_k^{a_k}=z_{k-1}\}}$$ allows us to verify easily the following: if
$h={\displaystyle {m\over a_1a_2...a_n}}$ is a character of the
group $\Sigma_{\boldmath{a}}$, then $(z, h)=z_{n+1}^m, \ z=(z_1,
z_2, \dots, z_n, \dots)\in G$.

\section{Main results}

Let $\boldmath{a}=(a_1,a_2,...)$ be a sequence of integers where all
$a_j>1$. Consider the group $ X=H_{\boldmath{a}}.$ Note that numbers
$$ e_{\pm 0}=\pm 1, \quad e_{\pm 1}=\pm {1\over a_1}, \quad
e_{\pm 2}\pm {1\over a_1a_2},\quad ..., \quad e_{\pm n}=\pm {1\over
a_1...a_n}, \quad... $$ are natural generators of the group $X$.

\noindent We consider on $X$ a distribution $\mu$ of the form

\begin{equation}\label{m1}
    \mu\{e_{\pm j}\}= {q_j\over 2},\quad \sum_{j=0}^{\infty} q_j =1,\quad q_j\geq 0.
\end{equation}
The distribution $\mu$ defines a random walk on $X$.

For a compact group $Y$ we suppose that a Haar measure $m_Y$ is
normalized in such a way that $m_Y(Y)=1$.

In the following theorem we obtain sufficient conditions for
recurrence of a random walk defined by a distribution $\mu$ on the
group $X=H_{\boldmath{a}}$.

\bigskip

\textbf{Theorem 1.} \textit{Let $X=H_{\boldmath{a}}$, where
$\boldmath{a}=(a_1,a_2,...), a_j>1$. Let $\mu$ be a distribution on
$X$ of the form $(\ref{m1})$. Consider on $X$ a random walk defined
by $\mu$. The condition
\begin{equation}\label{t1.6}
    \sum_{n=1}^{\infty} {1\over a_1...a_{n} \sqrt{q_{n}+
 q_{n+1}+...}}=\infty
\end{equation}
is sufficient for recurrence of a random walk defined by $\mu$ on
$X$.}

\bigskip

\textbf{ Proof.} The proof of Theorem 1 is based on Theorem A. The character group of the group $X$ is topologically isomorphic to the group $\Sigma_{\boldmath{a}}$. In order not to complicate the notation, we will assume that $Y=\Sigma_{\boldmath{a}}$. It is convenient to consider the realization of the $\boldmath{a}$-adic solenoid as a subgroup in $\mathbb{T}^{\aleph_0}$. Then each element of $Y$ is a sequence $y=(y_1,y_2,...)$, where $y_n\in
\mathbb{T}$, $y_n^{a_n}=y_{n-1}$. Put $y_n=e^{ib_n}$. Thus the sequence $(y_1,y_2,...)$ corresponds to a sequence $(b_1,b_2,...)$. Numbers $b_n$ are defined modulo
$2\pi$. We will take these numbers either in the interval $[0,2\pi)$ or in the interval $[-\pi,\pi)$ depending on how it is convenient for us. Regardless of an interval in which we take numbers
$b_n$, misunderstandings will not arise. Note that
\begin{equation}\label{t1.1.1.1}
    a_{n+1}b_{n+1}=b_n (mod \ 2\pi).
\end{equation}
Note also that $(e_{\pm j},y)= e^{\pm ib_{j+1}}=\cos b_{j+1}\pm
i \sin b_{j+1}$. Then
\begin{equation}\label{t1.1}
    \hat\mu(y)= \int_X (x,y) d \mu(x)={1\over 2}\sum_{j=0}^{\infty}
q_j (e_j,y) +{1\over 2}\sum_{j=0}^{\infty} q_j (e_{-j},y) =
\sum_{j=0}^{\infty} q_j \cos b_{j+1}.
\end{equation}

We will build by induction on the group $Y$ a system of non overlapping sets
$E_0, E_1, ..., E_n,...$ such that
\begin{equation}\label{t1.8}
    \sum_{n=0}^{\infty} m_Y(E_n)=1.
\end{equation}
Put
    $$ E_0=\left\{y\in Y:\ b_1 \not\in \left[-{\pi\over a_1}, {\pi\over a_1}\right] \right\}.$$
Using the invariance of the Haar measure we obtain that
\begin{equation}\label{t1.9}
    m_Y(E_0)={a_1-1\over a_1}.
\end{equation}
Note
    $$ Y\setminus E_0=\{y\in Y:\ b_1 \in \left[-{\pi\over a_1}, {\pi\over a_1}\right],
    b_2 \in \left[{2\pi k\over a_2}-{\pi\over a_1a_2}, {2\pi k\over a_2}+{\pi\over a_1a_2}\right] (k=0,1,...,a_2-1),...$$
    $$ ...,
    b_n \in \left[{2\pi k\over a_n}-{\pi\over a_1...a_n},
    {2\pi k\over a_n}+{\pi\over a_1...a_n}\right] (k=0,1,...,a_n-1), ... \} $$
See Figure 1 for the case when each $a_j=3$.

\begin{figure}
{\includegraphics[width=1\textwidth]{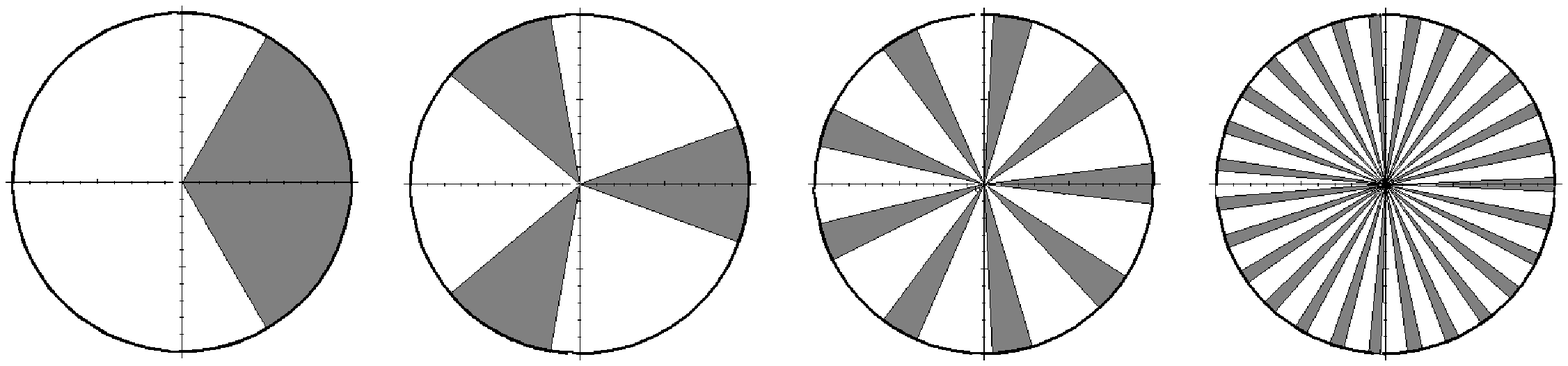}}
\caption{The domains of variation of elements $b_1$, $b_2$,
$b_3$, $b_4$ respectively in the subset $Y\setminus E_0$  for the case when each
$a_j=3$. }
\end{figure}

Put
\begin{equation}\label{t1.9}
    E_1=\left \{y\in Y\setminus E_0:\ b_2 \not\in \left[-{\pi\over a_1a_2}, {\pi\over
    a_1a_2}\right]\right \}=$$ $$
    \left\{y\in Y\setminus E_0:\ b_1 \in \left[-{\pi\over a_1}, {\pi\over a_1}\right],
    b_2 \in \left[{2\pi k\over a_2}-{\pi\over a_1a_2}, {2\pi k\over a_2}+{\pi\over a_1a_2}\right] (k=1,...,a_2-1)\right\}.
\end{equation}
It is easy to see that $m_Y(E_1)={a_2-1\over a_1a_2}$ (see Figure 2
for the case when each $a_j=3$).

\begin{figure}
{\includegraphics[width=1\textwidth]{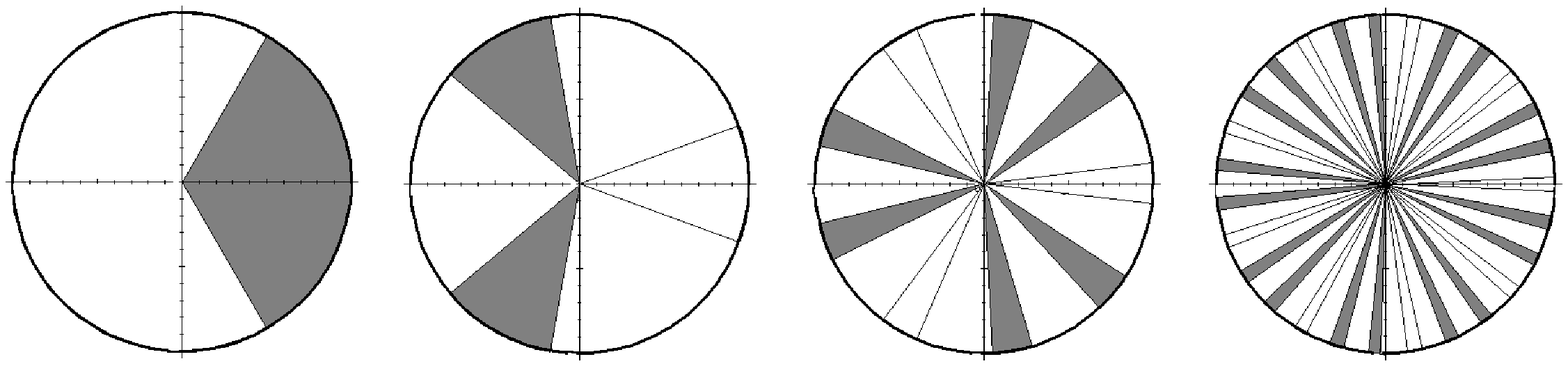}}
\caption{The domains of variation of elements $b_1$, $b_2$, $b_3$,
$b_4$ respectively in the subset $E_1$ for the case when each
$a_j=3$. }
\end{figure}

We define by induction a sequence of sets
\begin{equation}\label{t1.1.1}
    E_n=\left\{y\in Y\setminus
    \bigcup_{j=0}^{n-1}  E_j:\ b_{n+1} \not\in \left[-{\pi\over
    a_1...a_{n+1}}, {\pi\over a_1...a_{n+1}}\right]  \right\}.
\end{equation}

\noindent Actually sets $E_n$ can be defined by means of the
coordinate $b_{n+1}$. For a better understanding we note that

\begin{equation}\label{t1.1.1n}
    E_n=\left\{y\in Y: b_{n+1} \in \left[{2\pi k\over a_{n+1}}-{\pi\over a_1a_2...a_{n+1}}, {2\pi k\over a_{n+1}}+{\pi\over a_1a_2...a_{n+1}}\right] (k=1,...,a_{n+1}-1)\right\}=$$ $$=\{y\in Y\setminus \bigcup_{j=0}^{n-1}  E_j:\ b_1 \in \left[-{\pi\over a_1}, {\pi\over
    a_1}\right], b_2 \in \left[-{\pi\over a_1a_2}, {\pi\over
    a_1a_2}\right],..., b_n \in \left[-{\pi\over a_1a_2...a_n}, {\pi\over
    a_1a_2...a_n}\right],$$ $$
    b_{n+1} \in \left[{2\pi k\over a_{n+1}}-{\pi\over a_1a_2...a_{n+1}}, {2\pi k\over a_{n+1}}+{\pi\over a_1a_2...a_{n+1}}\right] (k=1,...,a_{n+1}-1)\}.
\end{equation}
See Figure 3 for the case when each $a_j=3$.

Using (\ref{t1.1.1}) and invariance of the Haar measure it is easy
to verify by induction that
\begin{equation}\label{t1.12}
    m_Y(E_n)={a_{n+1}-1\over a_1...a_{n+1}}
\end{equation}
for all $n=0,1,2,...$.

\begin{figure}
{\includegraphics[width=1\textwidth]{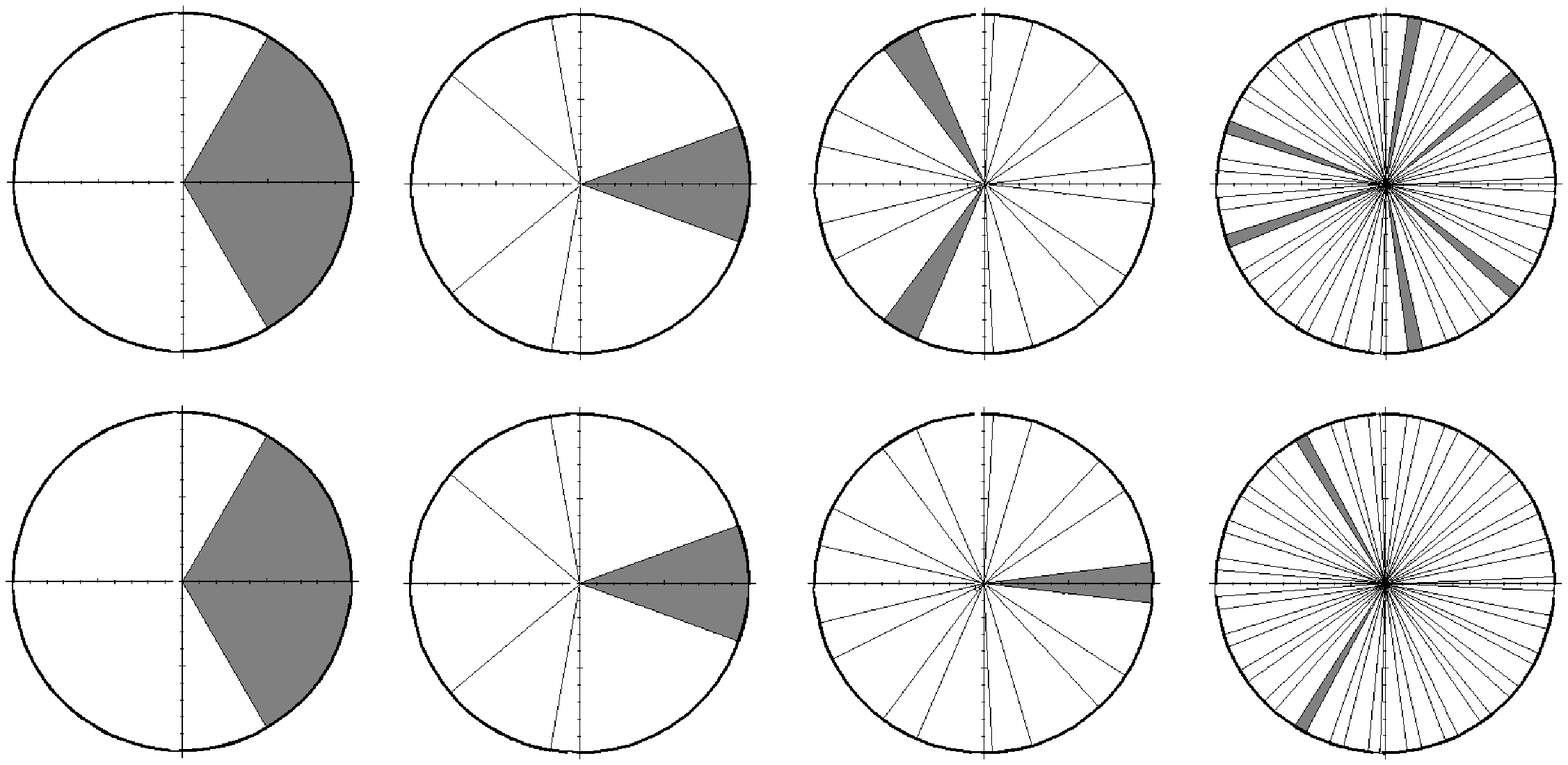}}
\caption{The domains of variation of elements $b_1$, $b_2$, $b_3$,
$b_4$ respectively in the subsets $E_2$ and $E_3$ for the case when
each $a_j=3$. }
\end{figure}

By construction, we have obtained that $E_i\cap E_j=\emptyset$ for
$i\neq j$. We have
$$\sum_{n=0}^{\infty} m_Y(E_n)={a_{1}-1\over a_1}+{a_{2}-1\over
a_1a_2}+{a_3-1\over a_1a_2a_3}+... = 1-{1\over a_1}+{1\over
a_1}-{1\over a_1a_2}+{1\over a_1a_2}-{1\over a_1a_2a_3}+...=1.$$

Let $\alpha_n$ be a number such that $0< \alpha_n <1$.

Put $$A_n=\left\{ y\in E_n:\  b_1 \in \left(-{\pi \alpha_n \over
a_1}, {\pi \alpha_n\over a_1}\right)\right\}=\left\{ y\in E_n:\  b_2
\in \left(-{\pi \alpha_n \over a_1a_2}, {\pi \alpha_n\over
a_1a_2}\right)\right\}=...$$
$$=\left\{ y\in E_n:\ b_k \in \left(-{\pi \alpha_n \over a_1a_2...a_k}, {\pi
\alpha_n\over a_1a_2...a_k}\right)\right\}=...$$ Since the measure
of $E_n$ coincides with the measure of projections on the $(n+1)$
circle, it is obvious that
\begin{equation}\label{t1.11}
    m_Y(A_n)=\alpha_n m_Y(E_n)= \alpha_n {a_{n+1}-1\over
a_1...a_{n+1}}.
\end{equation}

We evaluate from below the sum of the following series
\begin{equation}\label{t1.7}
    \sum_{n=0}^{\infty}
\int_{E_n} {1\over 1-\hat\mu(y)} dm_Y(y).
\end{equation}
We have
\begin{equation}\label{t1.2}
    \sum_{n=0}^{\infty} \int_{E_n} {1\over 1-\hat\mu(y)} dm_Y(y) \geq
\sum_{n=1}^{\infty} \int_{A_n} {1\over 1-\hat\mu(y)} dm_Y(y).
\end{equation}

We evaluate from above $1-\hat\mu(y)$ for $y\in A_n$. We have
$$1-\hat\mu(y)=1 - \sum_{j=0}^{\infty} q_j \cos b_{j+1}
=\sum_{j=0}^{\infty} q_j(1- \cos b_{j+1}) \leq $$

\begin{equation}\label{t1.2.1} \leq q_0 \left(1-\cos {\pi \alpha_n \over a_1}\right) + q_1 \left(1-\cos
{\pi \alpha_n \over a_1a_2}\right)+...+ q_{n-1} \left(1-\cos {\pi
\alpha_n \over a_1...a_{n}}\right)+  2\left( q_n+ q_{n+1}+...\right)
\end{equation} Since

\begin{equation}\label{t1.3}
    {2t^2\over \pi^2}\leq 1-\cos t\leq {t^2\over 2}, \quad t\in [0,{\pi\over 2}],
\end{equation}
we can continue evaluation (\ref{t1.2.1}) in the following way.

$$1-\hat\mu(y) \leq
{\pi^2\alpha_n^2 \over a_1^2}  \left(1+{1 \over a_2^2} + {1 \over
a_2^2 a_3^2} + ... + {1 \over a_2^2...a_n^2}\right) + 2\left( q_n+
q_{n+1}+...\right)
$$
\begin{equation}\label{t1.3.1}
\leq C_1 (\alpha_n^2 + q_n+ q_{n+1}+...),
\end{equation}
where $C_1$ is a constant which does not depend on $n$.

Taking into account (\ref{t1.3.1}) and (\ref{t1.11}), we can
continue evaluation (\ref{t1.2}):
\begin{equation}\label{t1.4}
    \sum_{n=0}^{\infty} \int_{E_n} {1\over 1-\hat\mu(y)} dm_Y(y) \geq
\sum_{n=1}^{\infty} \int_{A_n} {1\over 1-\hat\mu(y)} dm_Y(y)\geq C_2
\sum_{n=0}^{\infty}  {\alpha_n (a_{n+1}-1)\over a_1...a_{n+1}
(\alpha_n^2 + q_n+
 q_{n+1}+...)}$$
 $$ \geq C_3
\sum_{n=0}^{\infty}  {\alpha_n \over a_1...a_{n} (\alpha_n^2 + q_n+
 q_{n+1}+...)},
\end{equation}
where $C_2, C_3$ are constants which do not depend on $n$.

Note that for $a\in (0,1)$ the function $f(x)={x\over x^2+a}$ has a
maximum value on $[0,1]$ at $x=\sqrt{a}$. Taking this into account
from (\ref{t1.4}) we find that
\begin{equation}\label{t1.5}
    \sum_{n=0}^{\infty} \int_{E_n} {1\over 1-\hat\mu(y)} dm_Y(y) \geq
    C_4
\sum_{n=0}^{\infty}  {1 \over a_1...a_{n} \sqrt{q_n+
 q_{n+1}+...}},
\end{equation}
where $C_4$ is a constant which does not depend on $n$.

Now we can prove that condition (\ref{t1.6}) is sufficient for
recurrence of a random walk defined by the distribution $\mu$.
Indeed, suppose that condition (\ref{t1.6}) is satisfied and a
random walk defined by the distribution $\mu$ is transient. Then
Theorem A implies
$$\int_Y {1\over 1-\hat\mu(y)} dm(y)<\infty.$$ But then
$$ \int_Y {1\over 1-\hat\mu(y)} dm(y)=
\sum_{n=0}^{\infty} \int_{E_n} {1\over 1-\hat\mu(y)} dm(y)<\infty.
$$
It follows from (\ref{t1.5}) that series in (\ref{t1.6}) converges.
$\blacksquare$

\bigskip

\textbf{Remark 1.} If we put in Theorem 1 all $a_j=p$ where $p$  is
prime, we obtain Theorem 3 of the article \cite{Ferage-Molchanov}
for a random walk defined on the group $H_p$.

\bigskip

Thus, we have a sufficient condition for the recurrence of a random
walk defined on a subgroup of the group of rational numbers
$\mathbb{Q}$, not isomorphic to $\mathbb{Z}$. In the following
theorem we obtain a necessary condition for the recurrence of a
random walk defined on a subgroup of a group of rational numbers
$\mathbb{Q}$.

We need the following well-known property of characteristic
functions (see e.g. \cite[\S2]{Fe2008}).

\bigskip

\textbf{Lemma 1.} \textit{Let $X$ be a second countable locally
compact  Abelian group, let $\mu$ be a distribution on $X$. The
following conditions are equivalent.}

\textit{(i) The support of the distribution $\mu$ is not contained
in any coset of some subgroup in $X$.}

\textit{(ii) $\{y\in Y: |\hat\mu(y)|=1\}=\{0\}$.}

\bigskip

\textbf{Theorem 2.} \textit{Let $X=H_{\boldmath{a}}$, where
$\boldmath{a}=(a_1,a_2,...), a_j>1$. Let $\mu$ be a distribution on
$X$ of the form $(\ref{m1})$. We suppose that in $(\ref{m1})$ all
$q_j>0$. The condition
\begin{equation}\label{t2.1}
    \sum_{n=1}^{\infty} {a_{n+1}\over a_1...a_{n} \sqrt{q_{n}}}=\infty
\end{equation}
is necessary for recurrence of a random walk defined by $\mu$ on
$X$.}

\bigskip

\textbf{Proof.} Let sets $E_n$ be such as in Theorem 1. We evaluate
from above the sum of seriis (\ref{t1.7}).

Note that the support of the distribution $\mu$ is not contained in
any coset of some subgroup in $X$. Then Lemma 1 implies that

\begin{equation}\label{m1.1}
    \{y\in Y: |\hat\mu(y)|<1\}=\{0\}.
\end{equation}
Therefore
    $$ {1\over 1-\hat\mu(y)}= \sum_{k=0}^{\infty} \hat\mu^k(y), \quad y\neq 0.$$
By the Lebesgue-Levi theorem we have
\begin{equation}\label{3.1}
    \sum_{n=0}^{\infty} \int_{E_n} {1\over 1-\hat\mu(y)} dm_Y(y)=
    \sum_{n=0}^{\infty} \int_{E_n} \sum_{k=0}^{\infty} \hat\mu(y)^k
    dm_Y(y)=$$ $$
    \sum_{n=0}^{\infty} \sum_{k=0}^{\infty} \int_{E_n}  \hat\mu(y)^k dm_Y(y).
\end{equation}

It follows from (\ref{t1.1}) that $\hat\mu(y)\leq q_0 \cos {\pi\over
a_1}+q_1+q_2+...=1-q_0(1-\cos {\pi\over a_1})$ for $y\in E_0$. Hence
$$ \int_{E_0} {1\over 1-\hat\mu(y)} dm(y) <
{a_1-1\over a_1 q_0(1-\cos {\pi\over a_1})}.$$ Besides taking into
account (\ref{t1.12}), we note that
$$ \sum_{n=1}^{\infty} \int_{E_n}  dm(y) ={1\over a_1}. $$

Hence
\begin{equation}\label{3.1.1}
    \sum_{n=0}^{\infty} \int_{E_n} \sum_{k=0}^{\infty} \hat\mu(y)^k
    dm_Y(y)=C_1+\sum_{n=1}^{\infty} \sum_{k=1}^{\infty} \int_{E_n}
    \hat\mu(y)^k dm_Y(y),
\end{equation}
where $C_1$ is a constant which does not depend on $n$.

We evaluate from above $\hat\mu(y)$ for $y\in E_n$. It is obvious
that

\begin{equation}\label{t2.1}
    \hat\mu(y)=\sum_{j=0}^{\infty} q_j \cos b_{j+1} \leq \sum_{j=0}^{n} q_j
\cos b_{j+1} + q_{n+1}+q_{n+2}+ ... $$
$$= \sum_{j=0}^{n} q_j \cos b_{j+1} + 1-q_{0}-...-q_{n}=
1-\sum_{j=0}^{n} q_j (1-\cos b_{j+1}).
\end{equation}

Since $y\in E_n$, it follows from (\ref{t1.1.1}) that $b_{n+1} \in
\left[{2\pi k\over a_{n+1}}-{\pi\over a_1...a_{n+1}}, {2\pi k\over
a_{n+1}}+{\pi\over a_1...a_{n+1}}\right] (k=1,...,a_{n+1}-1)$. Hence
for any $n$ the inequality
\begin{equation}\label{t2.1n}
    1-\cos b_{n+1}\geq 1-\cos
    \left({2\pi \over a_{n+1}}-{\pi\over a_1...a_{n+1}}\right)
\end{equation}
is fulfilled. If $a_{n+1}=2$ or $a_{n+1}=3$ than ${\pi\over 2}<
{2\pi \over a_{n+1}}-{\pi\over a_1...a_{n+1}}< {\pi}$. Hence, taking
into account (\ref{t1.3}) and (\ref{t2.1n}), we have
\begin{equation}\label{t2.1n.1}
    1-\cos b_{n+1}\geq 1> {2\over a_{n+1}^2}.
\end{equation}
If $a_{n+1}>3$ then $0<{2\pi \over a_{n+1}}-{\pi\over
a_1...a_{n+1}}< {\pi\over 2}$. Hence, taking into account
(\ref{t2.1n}) and (\ref{t1.3}), we have
\begin{equation}\label{t2.1n.2}
    1-\cos b_{n+1}\geq {2\over \pi^2} \left( {2\pi\over a_{n+1}}-{\pi\over a_1...a_{n+1}} \right)^2\geq {2\over a_{n+1}^2}.
\end{equation}

Recall also that $b_n\in \left[-{\pi\over a_1...a_{n}}, {\pi\over
a_1...a_{n}}\right]$ for $y\in E_n$. Putting $t=b_{n}$ and taking
into account (\ref{t1.1.1.1}), we can rewrite estimate (\ref{t2.1})
in the following way

\begin{equation}\label{3.2}
\hat\mu(y) \leq 1-q_0(1-\cos a_{2}...a_{n} t)-q_1(1-\cos
a_{3}...a_{n} t)-...-q_{n-2}(1-\cos a_{n} t)-q_{n-1}(1-\cos t)-C_2
{q_n\over a_{n+1}^2},
\end{equation}
where $t\in \left[-{\pi\over a_1...a_{n}}, {\pi\over
a_1...a_{n}}\right]$.

Put $\tilde{t}=a_{2}...a_{n} t$. Then $\tilde{t}\in \left[-{\pi\over
a_1}, {\pi\over a_1}\right]$. We can rewrite estimate (\ref{3.2}) in
the following way
\begin{equation}\label{3.2.1}
\hat\mu(y) \leq 1-q_0(1-\cos \tilde{t})-q_1\left(1-\cos
{\tilde{t}\over a_2} \right)-...-q_{n-2}\left(1-\cos {\tilde{t}\over
a_2...a_{n-1}}\right)-q_{n-1}\left(1-\cos {\tilde{t}\over
a_2...a_{n}}\right)-C_2 {q_n\over a_{n+1}^2},
\end{equation}
where $\tilde{t}\in \left[-{\pi\over a_1}, {\pi\over a_1}\right]$.

Taking into account (\ref{t1.3}) we obtain
    $$ 1-\cos {\tilde{t}\over a_2...a_{j}}
    \geq {2\tilde{t}^2\over \pi^2 a_{2}^2...a_{j}^2}.$$

Thus we can continue inequality (\ref{3.2.1})

\begin{equation}\label{3.3}
\hat\mu(y) \leq 1-{2\tilde{t}^2\over \pi^2} \left(q_0 + {q_1\over
a_2^2}+{q_2\over a_2^2a_3^2}+ ... + {q_{n-1}\over a_2^2
a_3^2...a_n^2} \right) -C_2 {q_n\over a_{n+1}^2} \leq
e^{-C_3\tilde{t}^2-C_4 {q_n\over a_{n+1}^2}},
\end{equation}
where constants $C_3$ and $C_4$ are positive and do not depend on
$n$.

Note that if $y=(y_1,y_2,...)\in Y$ then $\tilde{t}=a_{2}...a_{n}
arg y_n$. Hence we can rewrite estimate (\ref{3.3}) in the folowing
way.
\begin{equation}\label{3.4}
\hat\mu(y) \leq e^{-C_3({a_{2}...a_{n}\arg y_n})^2-C_4 {q_n\over
a_{n+1}^2}}.
\end{equation}

Fix $n$. The mapping $f: Y \longrightarrow \mathbb{T}$ is defined by
the formula $f(y)=arg y_n$. It is easy to see that this mapping
transform the Haar measure on $Y$ into the Haar measure on
$\mathbb{T}$. We use the formula of the change of variable in the
integral. Then
\begin{equation}\label{3.5}
\int_{E_n}  \hat\mu^k(y) dm_Y(y) \leq \int_{E_n}
e^{-C_3k({a_{2}...a_{n}\arg y_n})^2-C_4 k {q_n\over a_{n+1}^2}}
dm_Y(y) = \int_{-\pi/a_1...a_n}^{\pi/a_1...a_n}
e^{-C_3k({{a_{2}...a_{n} t}})^2-C_4 k {q_n\over a_{n+1}^2}} dt.
\end{equation}

Next, we change variable $s=a_2...a_n t$. We have
\begin{equation}\label{3.6}
\int_{-\pi/a_1...a_n}^{\pi/a_1...a_n} e^{-C_3k({a_2...a_n t})^2-C_4
k {q_n\over a_{n+1}^2}} dt = {1\over
a_2...a_n}\int_{-\pi/a_1}^{\pi/a_1} e^{-C_3k s^2-C_4 k {q_n\over
a_{n+1}^2}} ds \leq C_5 {1\over a_1...a_n}  {e^{-C_4 k {q_n\over
a_{n+1}^2}}\over \sqrt{k}},
\end{equation}
where $C_5$ is a constant which does not depend on $n$ and $k$.

Thus
\begin{equation}\label{3.7}
    \sum_{n=1}^{\infty} \sum_{k=1}^{\infty} \int_{E_n}  \hat\mu(y)^k
    dm_Y(y)\leq C_5 \sum_{n=1}^{\infty} {1\over a_1...a_n}
    \sum_{k=1}^{\infty} {e^{-C_4 k {q_n\over a_{n+1}^2}}\over \sqrt{k}}.
\end{equation}

Since $\sum_{k=1}^{\infty} {e^{-a k}\over \sqrt{k}}=O({1\over
\sqrt{a}})$ (see e.g. \cite{Ferage-Molchanov}), we can continue
inequality (\ref{3.7}) in the following way

\begin{equation}\label{3.8}
    \sum_{n=1}^{\infty} \sum_{k=1}^{\infty} \int_{E_n}  \hat\mu(y)^k
    dm_Y(y)\leq C_6 \sum_{n=1}^{\infty} {a_{n+1}\over a_1...a_n
    \sqrt{q_n}},
\end{equation}
where $C_6$ is a constant which does not depend on $n$ and $k$.

Arguing in the same way as at the end of proof of Theorem 1, we
obtain that condition (\ref{t2.1}) is necessary for reccurence of a
random walk defined by $\mu$. $\blacksquare$

\bigskip

\textbf{Remark 2.} If we put in Theorem 2 all $a_j=p$ where $p$ is
prime, we obtain Theorem 4 of the article \cite{Ferage-Molchanov}
for a random walk defined on the group $H_p$.

\end{document}